\documentclass[11pt]{article}
\usepackage{amssymb,amsmath,latexsym}
\usepackage{color}

\usepackage{epsfig}

\newtheorem{thm}{Theorem}[section]

\newtheorem{corollary}[thm]{Corollary}
\newtheorem{lemma}[thm]{Lemma}
\newtheorem{remark}[thm]{Remark}

\newcommand{\pf}{\noindent{\bf Proof.} }
\def\qed{{\hfill $\Box$ \bigskip}}

\def\Xint#1{\mathchoice
{\XXint\displaystyle\textstyle{#1}}%
{\XXint\textstyle\scriptstyle{#1}}%
{\XXint\scriptstyle\scriptscriptstyle{#1}}%
{\XXint\scriptscriptstyle\scriptscriptstyle{#1}}%
\!\int}
\def\XXint#1#2#3{{\setbox0=\hbox{$#1{#2#3}{\int}$}
\vcenter{\hbox{$#2#3$}}\kern-.5\wd0}}

\def\dashint{\Xint-}

\newcommand\cbrk{\text{$]$\kern-.15em$]$}}
\newcommand\opar{\text{\,\raise.2ex\hbox{${\scriptstyle
|}$}\kern-.34em$($}}
\newcommand\cpar{\text{$)$\kern-.34em\raise.2ex\hbox{${\scriptstyle |}$}}\,}

\def\<{\langle}
\def\>{\rangle}

\newcommand\bR{\mathbb{R}}

\newcommand\bM{\mathbb{M}}
\newcommand\bE{\mathbb{E}}

\newcommand\cA{\mathcal{A}}
\newcommand\cB{\mathcal{B}}

\newcommand\cF{\mathcal{F}}
\newcommand\cG{\mathcal{G}}

\def\Re{\text{Re}}
\def\Im{\text{Im}}

\newcommand{\mysection}[1]{\section{#1}
\setcounter{equation}{0}}

\begin{document}

\title{\bf   A generalization of the Littlewood-Paley   inequality for the fractional Laplacian $(-\Delta)^{\alpha/2}$}

\date{}

\author{Ildoo Kim \qquad \hbox{\rm and} \qquad Kyeong-Hun Kim}

\author{Ildoo Kim \footnote{Department of Mathematics, Korea
University, 1 Anam-dong, Sungbuk-gu, Seoul, South Korea 136-701,
\,\, waldoo@korea.ac.kr.}\quad \hbox{\rm and} \quad Kyeong-Hun
Kim\footnote{Department of Mathematics, Korea University, 1
Anam-dong, Sungbuk-gu, Seoul, South Korea 136-701, \,\,
kyeonghun@korea.ac.kr. }}

\maketitle

\begin{abstract}
We prove a  parabolic version of the Littlewood-Paley  inequality
for the fractional Laplacian $(-\Delta)^{\alpha/2}$, where
$\alpha\in (0,2)$.

\vspace*{.125in}


\vspace*{.125in}


\end{abstract}

\section{Introduction}

Let $T_{2,t}$ be the semigroup corresponding to the heat equation
$u_t=\Delta u$ (see (\ref{06.09.1})). The classical Littlewood-Paley
inequality says for any $p\in (1,\infty)$ and $f\in L_p(\bR^d)$,
\begin{equation}
              \label{LP}
\int_{\bR^d} \left(\int^{\infty}_0 |\nabla T_{2,t}
f|^2dt\right)^{p/2} dx\leq N(p)\|f\|^p_p.
\end{equation}
 In \cite{kr94} and \cite{Kr01} Krylov
extended (\ref{LP}) by proving the following parabolic version in
which $H$ is a Hilbert space.
\begin{thm}
                    \label{krylov}
Let $H$ be a Hilbert space, $p\in [2,\infty), -\infty\leq a<b\leq
\infty$, $f\in L_p((a,b)\times \bR^d,H)$. Then
\begin{equation}
                      \label{eqn krylov}
 \int_{\bR^d}\int^b_a\left(\int^t_a|\nabla T_{2,t-s}f|^2_{H}\,ds\right)^{p/2}\,dtdx\leq
 N(p)
 \int_{\bR^d}\int^b_a |f|^p_{H}\,dtdx.
\end{equation}
\end{thm}

Let $\alpha\in (0,2)$. The main goal of this article is to  prove
(\ref{eqn krylov}) with $\partial^{\alpha/2}_x$ and $T_{\alpha,t}$
in place of $\nabla$ and $T_{2,t}$ respectively, where
$T_{\alpha,t}$ is the semigroup corresponding to the equation
$u_t=-(-\Delta)^{\alpha/2}u$. That is, we prove
\begin{thm}
                                             \label{1.1}
Let $H$ be a Hilbert space, $p\in [2,\infty), -\infty\leq a<b\leq
\infty$,  and $f$ be an $H$-valued function of $(t,x)$, then
\begin{equation}
                   \label{999}
\int_{\bR^d}\int^b_a\left
[\int^t_a|\partial^{\alpha/2}_xT_{\alpha,t-s}f(s,\cdot)(x)|^2_{H}
ds\right]^{p/2}dxdt\leq N(\alpha,p)\int_{\bR^d}\int^b_a
|f|^p_{H}\,dtdx.
\end{equation}
\end{thm}
 If $f(t,x)=f(x)$, then  (\ref{999}) easily leads to the Littlewood-Paley inequality
 (\ref{LP}) with $\partial^{\alpha/2}_x$ and $T_{\alpha,t}$ in place of $\nabla$ and $T_{2,t}$
  (see Remark \ref{remark 4}).

 Our motivation is as follows. For several decades,  the fractional Laplacian and partial
differential equations with the fractional Laplacian have been
studied by many authors, see for instance \cite{CS} and
\cite{Ste70}. Motivated by this, we were tempted to construct an
$L_p$-theory of stochastic partial differential equations of the
type
\begin{equation}
                 \label{eqn 0}
 du=-(-\Delta)^{\alpha/2} u \,dt +\sum_{k=1}^{\infty}f^kdw^k_t,
\quad u(0,x)=0.
\end{equation}
 Here
 $f=(f^1,f^2,\cdots)$ is an
$\ell_2$-valued random function of $(t,x)$, and $w^k_t$ are
independent one-dimensional Wiener processes.  It turns out that if
$f=(f^1,f^2,\cdots)$ satisfies certain measurability condition, the
solution of this problem is given by
\begin{equation}
                                   \label{eqn 222}
u(t,x)=\sum_{k=1}^{\infty}\int^t_0
T_{\alpha,t-s}f^k(s,\cdot)(x) dw^k_s,
\end{equation}
and by Burkholder-Davis-Gundy inequality (see \cite{kr95}), we have
\begin{equation}
                    \label{eqn 333}
\bE \int^T_0\|\partial^{\alpha/2}_xu(t,\cdot)\|^p_{L_p}dt \leq
N(p)\, \bE
\int^T_0\int_{\bR^d}\left[\int^t_0|\partial^{\alpha/2}_xT_{\alpha,t-s}f(s,\cdot)(x)|^2_{\ell_2}
ds\right]^{p/2}dxdt.
\end{equation}
Actually if $f$ is not random, then the reverse inequality also
holds. Thus to prove $\partial^{\alpha/2}_xu\in L_p$ and to get a
legitimate start of the $L_p$-theory of SPDEs of type (\ref{eqn 0}),
one has to estimate the right-hand side of (\ref{eqn 333}). Later,
we will see that (\ref{999}) implies that for any solution $u$ of
equation (\ref{eqn 0}),
\begin{equation}
                    \label{for fun}
\bE \int^T_0\|u(t,\cdot)\|^p_{H^{\alpha/2}_p}dt \leq
N(\alpha,p,T)\bE\int^T_0 \||f|_{\ell_2}\|^p_{L_p}ds,
\end{equation}
where $\|u\|_{H^{\alpha/2}_p}:=\|(1-\Delta)^{\alpha/4}u\|_{L_p}$.

As usual $\bR^{d}$ stands for the Euclidean space of points
$x=(x^{1},...,x^{d})$, $B_r(x) := \{ y\in \bR^d : |x-y| < r\}$  and
$B_r
 :=B_r(0)$.
 For  $\beta \in (0,1)$,
 and functions $u(x)$ we set
$$
\nabla_x u=(\frac{\partial}{\partial
x^1}u,\cdots,\frac{\partial}{\partial x^d}u), \quad
\partial_x^{\beta}u(x)=\cF^{-1}(|\xi|^{\beta}\hat{u}(\xi))(x)
$$
where
$\cF(f)(\xi)=\hat{f}(\xi):=\frac{1}{(2\pi)^d}\int_{\bR^d}e^{-i\xi\cdot
x}f(x)dx$ is the Fourier transform of $f$.  If we write $N=N(...)$,
this means that the constant $N$ depends only on what are in
parenthesis.

\section{Main Result}

In this section we introduce a  slightly extended version of Theorem
\ref{1.1}. Fix $\alpha \in (0,2)$ and let $p_{\alpha}(t,x)=p(t,x)$,
where $t>0$, denote the Fourier inverse transform of
$e^{-(2\pi)^\alpha t|\xi|^\alpha}$, that is,
$$
p(t,x):= \int_{\bR^d} e^{i \xi \cdot x }e^{-(2\pi)^\alpha
t|\xi|^\alpha} d\xi
$$
and $p(x) := p(1,x)$. For a suitable function $h$ and $t>0$,  define
\begin{equation}
                            \label{06.09.1}
T_th(x) := (p(t,\cdot) * h(\cdot))(x):=\int_{\bR^d} p(t,x-y)h(y)dy,
\end{equation}
$$
({-\Delta})^{\frac{\beta}{2}}h(x) :=\partial_x^{\beta} h:=
\cF^{-1}(|\xi|^{\beta} \cF(h)(\xi))(x).
$$
Then, for $\beta > 0$,
\begin{eqnarray}
\partial^{\beta}_x T_th(x)
&=&\cF^{-1}( |\xi|^{\beta} e^{-(2\pi)^\alpha t|\xi|^\alpha} \hat{h}(\xi))\nonumber\\
&=&\int_{\bR^d} e^{i\xi \cdot x}
|\xi|^{\beta} e^{-(2\pi)^\alpha  |t^{1/\alpha}\xi|^{ \alpha}} d\xi \ast h(x)\nonumber\\
&=&t^{-d/\alpha}\int_{\bR^d} e^{i \xi  \cdot t^{-1/{\alpha}} x} |t^{-1/\alpha}\xi|^{\beta}
e^{-(2\pi)^\alpha  |\xi|^\alpha} d\xi \ast h(x)\nonumber\\
&=&t^{-\beta/\alpha}\cdot t^{-d/\alpha} \phi_\beta(x /t^{1/\alpha})
\ast h(x) \label{eqn 06.12.1},
\end{eqnarray}
where
$$
 \phi_\beta(x) := \int_{\bR^d} |\xi|^{\beta} e^{i \xi \cdot x}
e^{-(2\pi)^\alpha |\xi|^{\alpha}} d\xi=
({-\Delta})^{\frac{\beta}{2}}p(x).
$$

The following  two lemmas are crucial in this article and are proved
in section \ref{section a priori}.

\begin{lemma}
                                               \label{lem3}
Denote $\hat{\phi}_\beta(\xi)=\cF(\phi_\beta)(\xi)$, then there
exists a constant $N=N(d,\alpha,\beta )>0$ such that
$$
|\hat{\phi}_\beta(\xi)| \leq N|\xi|^\beta, \quad |\xi|
|\hat{\phi}(\xi)| \leq N,
$$
$$
|\phi_\beta(x)|\leq N   \left( \frac{1}{|x|^{d+\beta}} \wedge 1
\right)\quad \text{and} \quad |\nabla \phi_\beta(x)| \leq N \left(
\frac{1}{|x|^{d+1+\beta}} \wedge 1 \right).
$$
\end{lemma}

\begin{lemma}
                                 \label{main lem}
 For each $\alpha\in (0,2)$ and $\beta > 0$, there exists a continuously differentiable function
  $\overline{\phi}_\beta(\rho)$ defined on $[0,\infty)$ such that for some positive constant
  $K$ which depends on $d,\alpha ,\beta$,
$$
|\phi_\beta(x)| + |\nabla \phi_\beta(x)| + |x||\nabla \phi_\beta(x)|
\leq \overline{\phi}_\beta(|x|), \quad \int_0^\infty
|\overline{\phi}_\beta'(\rho)|~d\rho \leq K,
$$
$$
\overline{\phi}_\beta(\infty)=0, \quad \int_r^\infty
|\overline{\phi}_\beta'(\rho)|\rho^d~d\rho \leq \frac{K}{r^{\beta}}
\quad \quad \forall r \geq (10)^{-1/\alpha}.
$$
\end{lemma}

To make our inequality slightly extended, we  consider convolutions
(see (\ref{eqn 06.12.1})) with  more general functions. Let
$\psi(x)$ be a $C^1(\bR^d)$  function such that $|\hat{\psi}(\xi)|
\leq K |\xi|^{\nu}$ for some $\nu >0 $, $|\xi|^\lambda
|\hat{\psi}(\xi)| \leq K$ for some $\lambda >0$, and assume that for
some $\delta \geq \frac{\alpha}{2}$, there exists a continuously
differentiable function $\overline{\psi}$ satisfying
$$
|\psi(x)|+|\nabla \psi(x)| + |x||\nabla \psi(x)| \leq
\overline{\psi}(|x|), \quad \int_0^\infty |\overline{\psi}'(\rho)|
~d\rho \leq K, \quad \overline{\psi}(\infty)=0
$$
and
\begin{equation}
                     \label{eqn 06.08}
\int_r^\infty |\overline{\psi}'(\rho)|\rho^d ~d\rho \leq
(K/r^\delta), \quad \forall r \geq (10)^{-1/\alpha}.
\end{equation}
By Lemma \ref{lem3} and Lemma \ref{main lem}, we know
$\phi_{\alpha/2}$ satisfies all the above assumptions. Define
$$
\Psi_th(x):=t^{-d/\alpha} \psi(\cdot /t^{1/\alpha}) \ast h(\cdot
)(x).
$$
For $f \in C_0^\infty (\bR^{d+1},H),~ t > a \geq -\infty$, and $x
\in \bR^d$, we define
$$
\cG_af(t,x):=[\int_a^t | \Psi_{t-s}f(s, \cdot)(x)|^2_H
~\frac{ds}{t-s}]^{1/2}, \quad \cG=\cG_{-\infty}.
$$

Here is our main result. The proof is given in  section \ref{proof
of theorem}.

\begin{thm}
                                           \label{main theorem}
Let $p \in [2,\infty)$, $-\infty \leq a < b \leq \infty$ and $f \in
C_0^\infty((a,b) \times \bR^d , H).
$
Then
\begin{eqnarray}
                                 \label{4022}
\int_{\bR^d} \int_a^b [\cG_a f(t,x)]^p ~dt dx \leq N
\int_{\bR^d} \int_a^b |f(t,x)|^p_H~dtdx,
\end{eqnarray}
where the constant $N$ depends only on
$d,p,\alpha,\nu,\lambda,\delta$ and $K$.
\end{thm}

\begin{remark}
                     \label{remark 22}
Take $\psi = \phi_{\alpha/2}$, $\nu=\delta=\alpha/2$, $\lambda=1$,
$a=0$ and $b=T$, then (\ref{4022}) implies
\begin{equation}
                         \label{eqn 6.10.5}
\int_{\bR^d} \int_0^T
[\int_0^t|\partial^{\alpha/2}_xT_{\alpha,t-s}f(s, \cdot)(x)|^2_H
ds]^{p/2} ~dt dx \leq N \int_{\bR^d} \int_0^T |f(t,x)|^p_H~dtdx.
\end{equation}
\end{remark}

\begin{remark}
                           \label{remark 4}
Note that  inequality (\ref{lem3}) with $\partial^{\alpha/2}_x$ and
$T_{\alpha,t}$ in place of $\nabla$ and $T_{2,t}$  is an easy
consequence of (\ref{eqn 6.10.5}). Indeed, take $T=2$ and
$f(t,x)=f(x)$. The left-hand side of (\ref{eqn 6.10.5}) is not less
than
$$
\int_{\bR^d}\int^2_1[\int^1_0|\partial^{\alpha/2}_xT_{\alpha,s}f(x)|^2_H\,ds]^{p/2}\,dtdx
=\int_{\bR^d}[\int^1_0|\partial^{\alpha/2}_xT_{\alpha,s}f(x)|^2_H\,ds]^{p/2}\,dx.
$$
Thus it follows that
$$
\int_{\bR^d}[\int^1_0|\partial^{\alpha/2}_xT_{\alpha,s}f(x)|^2_H\,ds]^{p/2}\,dx\leq
N\int_{\bR^d}\|f\|^p_{H}\,dx,
$$
and  the self-similarity
$(\partial^{\alpha/2}_xT_{\alpha,s}f(c\,\cdot))(x)=c^{\alpha/2}(\partial^{\alpha/2}_xT_{\alpha,c^{\alpha}s}f)(cx)$
allows one to replace the upper limit $1$ by infinity with the same
constant $N$.
\end{remark}

\section{Preliminary estimates on $({-\Delta})^{\beta/2}p(t,x)$}
                     \label{section a priori}

In this section we study the upper bound of
$|({-\Delta})^{\beta/2}p(t,x)|$ and
$|\nabla({-\Delta})^{\beta/2}p(t,x)|$, and then we prove Lemma
\ref{lem3} and Lemma \ref{main lem}.  Actually the arguments in this
section allow one to get the upper bound of
$|D^m({-\Delta})^{\beta/2}p(t,x)|$ for any $m\geq 0$.

\begin{lemma}
                                            \label{lem1}
 There exists a constant $N=N(d,\alpha ,\beta)>0$ such that
\begin{equation}
                         \label{eqn 2.1}
|({-\Delta})^{\frac{\beta}{2}}p(x)|
 \leq  \frac{N}{|x|^{d+\beta}}.
\end{equation}
\end{lemma}

\pf See \cite{fe} for $d=1$ and \cite{kol} for $d \geq 2$. Actually
in \cite{fe}, (\ref{eqn 2.1}) is given only for $\beta=0$. Also in
\cite{kol}, $(-\Delta)^{\frac{\beta}{2}}p(x)$ is estimated in terms
of power series (Proposition 2.2), however the series does not
converge if $\alpha >1$. For these reasons, we give a detailed
proof. Also some inequalities obtained  in this proof will be used
in the proof of Lemma \ref{lem2}.

For $d=1$, since $|\xi|$ is an even function, we have
\begin{eqnarray}
 ({-\Delta})^{\frac{\beta}{2}}p(x)
\notag &=&\int_{\bR} |\xi|^{\beta}e^{i \xi  x }e^{-(2\pi)^\alpha |\xi|^\alpha} d\xi\\
\notag &=&2 \Re \int_{0}^\infty \xi^{\beta}e^{i \xi  x }e^{-(2\pi)^\alpha \xi^\alpha} d\xi\\
\label{2.2}&=&2 \Re \frac{1}{x^{1+\beta}}\int_{0}^\infty \xi^{\beta}e^{i \xi }e^{-(2\pi)^\alpha (\xi/x)^\alpha} d\xi.
\end{eqnarray}
Assume  $0 < \alpha \leq 1$. Consider the integrand as a function of
the complex variable $\xi$. Since the integrand in (\ref{2.2}) is
analytic in the complement of the non-positive real half line and is
continuous at zero, if we take principal branch cut, for $N>0$, the
path integration is zero on the closed path
$$
\gamma_N(t) :=
\begin{cases}
t \quad & \text{if}~ 0 \leq t \leq N \\
N+i(t-N) \quad & \text{if} ~N \leq t \leq 2N \\
3N-t+iN \quad & \text{if} ~2N \leq t \leq 3N \\
i(4N-t) \quad & \text{if}~ 3N \leq t \leq 4N.
\end{cases}
$$
By letting $N \to \infty$, one can move the path of integration to
the positive imaginary axis, and gets (note $|e^{-(2\pi)^\alpha
(i\xi/x)^\alpha}| \leq 1$)
$$
|({-\Delta})^{\frac{\beta}{2}}p(x)| \label{2.31} = \left|2 \Re
\frac{1}{x^{1+\beta}}\int_{0}^\infty (i\xi)^{\beta}e^{- \xi }
e^{-(2\pi)^\alpha (i\xi/x)^\alpha} i\, d\xi\right|  \leq
\frac{2}{|x|^{1+\beta}}\int^{\infty}_0 \xi^{\beta}e^{-\xi}d\xi
 \leq  \frac{N}{|x|^{1+\beta}}.
$$
If $1 < \alpha  <2$, we use  another closed path
$$
\gamma_N(t) :=
\begin{cases}
 t \quad & \text{if}~ 0 \leq t \leq N \cos \frac{\pi}{2\alpha} \\
N \cos \frac{\pi}{2\alpha}
+i\sin\frac{\pi}{2\alpha}(\frac{t}{\cos \frac{\pi}{2\alpha}}-N)
\quad & \text{if} ~N\cos \frac{\pi}{2\alpha} \leq t \leq 2N\cos \frac{\pi}{2\alpha} \\
(3N-\frac{t}{\cos \frac{\pi}{2\alpha}})e^{i\frac{\pi}{2\alpha}}
 \quad & \text{if}~ 2N\cos \frac{\pi}{2\alpha} \leq t \leq 3N\cos \frac{\pi}{2\alpha}.
\end{cases}
$$
Thanks to the path integration on the above path, which looks like
formally replacing  $\xi$ by $\xi e^{i\frac{\pi}{2\alpha}}$, we get
(since $|e^{-(2\pi)^\alpha (\xi
e^{i\frac{\pi}{2\alpha}}/x)^\alpha}|= 1$)
\begin{eqnarray*}
|({-\Delta})^{\frac{\beta}{2}}p(x)|
&\leq& |2 \Re \frac{1}{x^{1+\beta}}\int_{0}^\infty
(\xi e^{i\frac{\pi}{2\alpha}})^{\beta}e^{ i\xi e^{i\frac{\pi}{2\alpha}} }
e^{-(2\pi)^\alpha (\xi e^{i\frac{\pi}{2\alpha}}/x)^\alpha} e^{i\frac{\pi}{2\alpha}}\,d\xi| \\
&\leq & \frac{2}{|x|^{1+\beta}}\int^{\infty}_0 \xi^{\beta}e^{-\xi
\sin \frac{\pi}{2\alpha}}\,d\xi \leq\frac{N}{|x|^{1+\beta}}.
\end{eqnarray*}

Next, let $d \geq 2$. Since the function $({-\Delta})^{\frac{\beta}{2}}p(x) $ is radial, we may assume $x= (|x|,\ldots,0)$, and if we denote the surface of the $d$-dimensional unit ball by $S^{d-1}$ and the surface measure by $d\sigma$, then from the spherical coordinate we have
\begin{eqnarray}
\notag ({-\Delta})^{\frac{\beta}{2}}p(x)
\notag &=&\int_{\bR^d} |\xi|^{\beta}e^{i \xi^1  |x| }e^{-(2\pi)^\alpha |\xi|^\alpha} d\xi\\
\notag &=&\int_{\bR^d} |\xi|^{\beta}\cos(  \xi^1  |x| )e^{-(2\pi)^\alpha |\xi|^\alpha} d\xi\\
\notag &=&\int_{0}^\infty r^{\beta+d-1}
\int_{S^{d-1}} \cos( r \sigma^1   |x| )e^{-(2\pi)^\alpha |r |^\alpha} d\sigma dr.
\end{eqnarray}
Furthermore we can express $\sigma \in S^{d-1}$ as $\sigma =
(\cos\theta, \phi \sin \theta)$ with $\theta \in [0,\pi]$ and $\phi
\in S^{d-2}$, and get
\begin{eqnarray*}
\notag ({-\Delta})^{\frac{\beta}{2}}p(x)
\notag &=&\int_{0}^\infty r^{\beta+d-1} \int_0^{\pi}  \sin^{d-2} (\theta) \int_{S^{d-2}} \cos( r \cos\theta  |x| )e^{-(2\pi)^\alpha |r |^\alpha} d\phi d\theta dr\\
\notag &=&A_{d-2}\int_{0}^\infty r^{\beta+d-1} \int_0^{\pi}  \sin^{d-2} (\theta) \cos( r \cos\theta   |x| ) e^{-(2\pi)^\alpha |r |^\alpha}  d\theta dr,
\end{eqnarray*}
where $A_{d-2}$ is the area of $S_{d-2}$ and $A_0:=1$. By the
changes of variables $r|x| \to r$ and  $t = \cos \theta$,
\begin{eqnarray}
({-\Delta})^{\frac{\beta}{2}}p(x) \notag &=&A_{d-2}
\frac{1}{|x|^{\beta+d}} \int_{0}^\infty r^{\beta+d-1} \int_0^{\pi}
\sin^{d-2} (\theta) \cos( r \cos\theta    ) e^{-(2\pi)^\alpha
(r/|x|)^\alpha}  d\theta dr \nonumber\\
  &=&A_{d-2} \frac{1}{|x|^{\beta+d}} \int_{0}^\infty r^{\beta+d-1}
  \int_{-1}^{1}   \cos( r t    ) e^{-(2\pi)^\alpha  (r/|x|)^\alpha}  (1-t^2)^{(d-3)/2} dt dr \label{2.3}.
\end{eqnarray}
To proceed further, we use  Bessel function $J_n(z)$ and  Whittaker
function $W_{0,n}(z)$.  For any complex $z$ that is not  negative
real and any real $n > -\frac{1}{2}$, define
$$
J_n(z) := \frac{(\frac{1}{2}z)^n}{\Gamma(n+\frac{1}{2})\sqrt{\pi}}
\int_{-1}^1 (1-t^2)^{n-1/2} \cos (zt) dt,
$$

\begin{equation}
                                               \label{5141}
W_{0,n}(z):= \frac{e^{-z/2}}{\Gamma(n+\frac{1}{2})}
\int_0^\infty [t(1+t/z)]^{n-1/2}e^{-t}dt
\end{equation}
where $\arg z$ is understood to take its principle value, that is,
$|\arg z| < \pi$. It is known (see, for instance, \cite{whit} p.346,
p.360 and \cite{wang} p.314) that  the two functions are related by
the formula
$$
J_n(z) = \frac{1}{\sqrt{2\pi z}} \left( \exp \{ \frac{1}{2}(n+
\frac{1}{2}) \pi i\} W_{0,n}(2 i z) +\exp \{ -\frac{1}{2}(n+
\frac{1}{2}) \pi i\} W_{0,n}(-2 i z) \right).
$$
In particular,  if $z$ is a positive  real number,
\begin{equation}
\label{2.4} J_n(z) = 2 \Re\left[ \frac{1}{\sqrt{2\pi z}} \exp \{
\frac{1}{2}(n+\frac{1}{2})\pi i \} W_{0,n}(2 i z) \right].
\end{equation}
We also know (see, for instance, \cite{whit} p. 343)
\begin{eqnarray}
\label{2.5} W_{0,n}(z) = e^{-\frac{1}{2}z}  \{ 1+ O(z^{-1})\}.
\end{eqnarray}
Due to (\ref{5141}) and (\ref{2.4}), from (\ref{2.3}) we have
\begin{eqnarray}
\notag && ({-\Delta})^{\frac{\beta}{2}}p(x)\\
\notag
&=& \frac{A_{d-2}}{|x|^{\beta+d}} \int_{0}^\infty r^{\beta+d-1}
 \int_{-1}^{1}   \cos( r t    ) e^{-(2\pi)^\alpha  (r/|x|)^\alpha}  (1-t^2)^{(d-3)/2} dt dr\\
 &=& \frac{A_{d-2}}{|x|^{\beta+d}} \int_{0}^\infty r^{\beta + d/2}
2^{d/2-1}\Gamma(\frac{1}{2}(d-1)) \sqrt{\pi}J_{(d/2)-1}(r)  e^{-(2\pi)^\alpha (r/|x|)^\alpha} dr\label{eqn 6.4.5}\\
 &=& \frac{N(d)}{|x|^{\beta+d}} \Re \int_{0}^\infty r^{\beta +
(d-1)/2}
 \exp \{ \frac{1}{2}(\frac{d}{2}-\frac{1}{2})\pi i \} W_{0,(d/2)-1}(2 i r)
  e^{-(2\pi)^\alpha (r/|x|)^\alpha} dr, \label{2.66}
\end{eqnarray}
where $N(d):=2^{(d-1)/2}A_{d-2}\Gamma(\frac{1}{2}(d-1))$. From
definition (\ref{5141}) one easily checks that the integrand in
(\ref{2.66}) is analytic in the complement of the non-positive real
half line and is continuous at zero.

Let $0<\alpha \leq 1$. Remembering (\ref{2.5}) and doing the path
integration on an appropriate closed path, as in the case $d=1$, we
can change the path of integration in (\ref{2.66}) from the positive
real half line to the negative imaginary half line. Taking this new
path of integration, that is to say, formally replacing $r$ by
$-ir$,
 one gets (note $|e^{-(2\pi)^\alpha (-ir/|x|)^\alpha}|\leq 1$)
\begin{eqnarray*}
\notag && |({-\Delta})^{\frac{\beta}{2}}p(x)|\\
\notag &=& \frac{N(d)}{|x|^{\beta+d}}
 \left|\Re \int_{0}^\infty (-ir)^{\beta + (d-1)/2}
  \exp \{ \frac{1}{2}(\frac{d}{2}-\frac{1}{2})\pi i \} W_{0,(d/2)-1}(2 r)
  e^{-(2\pi)^\alpha (-ir/|x|)^\alpha} i\, dr \right|\\
\notag &\leq& \frac{N}{ |x|^{\beta+d}}  \int_{0}^\infty r^{\beta +
(d-1)/2}
 W_{0,(d/2)-1}(2 r)   dr \leq \frac{N}{ |x|^{\beta+d}}.
\end{eqnarray*}

Let $1 <\alpha <2$. Then $|e^{-i re^{-i\frac{\pi}{2\alpha}}}|\leq
e^{\frac{-r}{2}}$, and thus
$$
|W_{0,(d/2)-1}(2 ire^{-i\frac{\pi}{2\alpha}})|\leq
\frac{e^{\frac{-r}{2}}}{\Gamma(d/2-1/2)}\int_0^\infty
\left|[t(1+t/(2 ire^{-i\frac{\pi}{2\alpha}}))]^{(d-3)/2}
e^{-t}\right|dt.
$$
Note that if $d\geq 3$ then
$$|1+t/(2
ire^{-i\frac{\pi}{2\alpha}})|^{(d-3)/2}\leq |1+t/r|^{(d-3)/2},
$$
and if $d=2$ then
$$|1+t/(2 ire^{-i\frac{\pi}{2\alpha}})|^{-1/2}\leq
(1+t\sin\frac{\pi}{2\alpha}/(2r))^{-1/2}\leq 2(1+t/r)^{-1/2}.
$$
It follows that for any $r>0$, we have $|W_{0,(d/2)-1}(2
ire^{-i\frac{\pi}{2\alpha}})|\leq
 2W_{0,(d/2)-1}(r)$.

We do the path integration on a different closed path and change the
path of integration in (\ref{2.66}) from the positive real half line
to the half line $\{r e^{-i\frac{\pi}{2\alpha}}: r > 0\}$. Taking
this new path of integration, that is to say, formally replacing $r$
by $re^{-i\frac{\pi}{2\alpha}}$, one gets (note $|e^{-(2\pi)^\alpha
 (re^{-i\frac{\pi}{2\alpha}}/|x|)^\alpha}|=1$)
\begin{eqnarray*}
|({-\Delta})^{\frac{\beta}{2}}p(x)| &\leq&
\frac{N}{|x|^{\beta+d}}\int_{0}^\infty \left|
(re^{-i\frac{\pi}{2\alpha}})^{\beta + (d-1)/2}  W_{0,(d/2)-1}(2
ire^{-i\frac{\pi}{2\alpha}})  e^{-(2\pi)^\alpha
 (re^{-i\frac{\pi}{2\alpha}}/|x|)^\alpha}\right|\,dr\\
&\leq& \frac{N}{ |x|^{\beta+d}}  \int_{0}^\infty r^{\beta +
\frac{d-1}{2}}
 W_{0,(d/2)-1}( r)  dr \leq \frac{N}{ |x|^{\beta+d}}.
\end{eqnarray*}
The lemma is proved. \qed

\begin{remark}
                \label{remark 6.4}
In the proof of Lemma \ref{lem1} (see (\ref{2.2}) and (\ref{eqn
6.4.5})) we proved that for any $\beta\geq 0$,
\begin{equation}
                      \label{eqn 6.4.3}
\left|\int^{\infty}_0 \xi^{\beta}e^{i\xi}e^{-(2\pi)^\alpha
(\xi/x)^\alpha} d\xi \right|<N,  \quad \text{when}\,\,d=1,
\end{equation}
\begin{equation}
                      \label{eqn 6.4.4}
\left|\int_{0}^\infty r^{\beta + d/2} J_{(d/2)-1}(r)
e^{-(2\pi)^\alpha (r/|x|)^\alpha} dr\right|<N, \quad \text{when}\,\,
d\geq 2,
\end{equation}
where $N=N(\alpha,\beta,d)>0$ is independent of $x$.
\end{remark}


\begin{remark}
                                   \label{remark 6.14}
Even though (\ref{eqn 2.1}) is enough for our need, we believe it is
not sharp. Actually it is known (see \cite{B}) that if $\beta=0$,
then
$$
p(t,x) \sim \left(\frac{t}{|x|^{d+\alpha}} \wedge
t^{-d/\alpha}\right).
$$
\end{remark}

\begin{lemma}
                                  \label{lem2}
 There exists a constant $N=N(d,\alpha , \beta)>0$ such that
\begin{equation}
                          \label{4021}
|\nabla ({-\Delta})^{\frac{\beta}{2}}p(x) | \leq
N(\frac{1}{|x|^{\beta+d+1}} \vee \frac{1}{|x|^{\beta+d+\alpha+1}}).
\end{equation}
\end{lemma}

\pf Let $d=1$. By (\ref{eqn 6.4.3}),
\begin{eqnarray*}
| \frac{d}{dx}({-\Delta})^{\frac{\beta}{2}}p(x)|
&=&|\int_{\bR} i  \xi |\xi|^{\beta }e^{i \xi  x }e^{-(2\pi)^\alpha |\xi|^\alpha} d\xi|\\
&\leq & \frac{1}{|x|^{\beta+2}}|\int_{\bR}   \xi |\xi|^{\beta }e^{i \xi   }e^{-(2\pi)^\alpha |\xi/x|^\alpha} d\xi|\\
&=&\frac{2}{|x|^{\beta+2}}| \Im \int_{0}^\infty  \xi^{1+\beta }e^{i
\xi   }e^{-(2\pi)^\alpha |\xi/x|^\alpha} d\xi|  \leq
\frac{N}{|x|^{\beta+2}}.
\end{eqnarray*}

Let $d\geq 2$. From (\ref{eqn 6.4.5}) and the inequality
$$
|\frac{\partial}{\partial x_i} ({-\Delta})^{\frac{\beta}{2}}p(x)|
=|\frac{\partial}{\partial |x|}
({-\Delta})^{\frac{\beta}{2}}p(x)\frac{\partial |x|}{\partial
x_i}|\leq |\frac{\partial}{\partial |x|}
({-\Delta})^{\frac{\beta}{2}}p(x)|
$$
it easily follows that
\begin{eqnarray*}
&&|({-\Delta})^{\frac{\beta}{2}}p(x)|\\
&\leq& \frac{N_1}{|x|^{\beta+d+1}} | \int_{0}^\infty (r)^{\beta +
d/2} J_{d/2-1}(r)e^{-(2\pi)^\alpha (r/|x|)^\alpha}
dr|\\
&&+\frac{N_2}{|x|^{\beta+d+\alpha+1}}| \int_{0}^\infty (r)^{\beta +
d/2+\alpha} J_{d/2-1}(r)e^{-(2\pi)^\alpha (r/|x|)^\alpha} dr|.
\end{eqnarray*}
Thus by (\ref{eqn 6.4.4}),
$$
|({-\Delta})^{\frac{\beta}{2}}p(x)|\leq N(\frac{1}{|x|^{\beta+d+1}}
\vee \frac{1}{|x|^{\beta+d+\alpha+1}}).
$$
The lemma is proved. \qed

({\bf{Proof of Lemma \ref{lem3}}})

  First two assertions come
from the fact
$$
\cF(\phi_\beta(x))(\xi) =\cF( \int_{\bR^d} |\eta|^{\beta} e^{i\eta
\cdot x} e^{-(2\pi )^\alpha |\eta|^{\alpha}}~d\eta)(\xi)
=|\xi|^{\beta} e^{-(2\pi )^\alpha |\xi|^{\alpha}}.
$$

Next, observe that
$$
|\phi_\beta(x)|=|(-\Delta)^{\frac{\beta}{2}}p(x)| =
|\int_{\bR^d}|\xi|^{\beta}e^{i\xi \cdot x} e^{-(2\pi)^\alpha
|\xi|^\alpha}d\xi|  \leq  \int_{\bR^d}|\xi|^{\beta}
e^{-(2\pi)^\alpha |\xi|^\alpha}d\xi<\infty.
$$
 Similarly,
$$
|\nabla \phi_\beta(x)|\leq \int_{\bR^d}|\xi|^{\beta+1}
e^{-(2\pi)^\alpha |\xi|^\alpha}d\xi<\infty.
$$
Therefore, by  Lemma \ref{lem1} and Lemma \ref{lem2}, there exists a
constant $N(d,\alpha,\beta)>0$ such that
\begin{eqnarray*}
&|\phi_\beta (x)|\leq N   \left( \frac{1}{|x|^{d+\beta}} \wedge 1
\right), \quad |\nabla \phi_\beta(x)| \leq N    \left(
\frac{1}{|x|^{d+1+\beta}} \wedge 1 \right).
\end{eqnarray*}
The lemma is proved. \qed

{\bf{(Proof of Lemma \ref{main lem})}}\\
By the inequalities in Lemma \ref{lem3}, we have
\begin{eqnarray*}
|\phi_\beta(x)| + |\nabla \phi_\beta(x)| + |x||\nabla \phi_\beta(x)|
\leq N \left(\frac{1}{|x|^{d+\beta}} \wedge 1 \right).
\end{eqnarray*}
Define
$$
\overline{\phi}_\beta(\rho)=\begin{cases} &\frac{N}{\rho^{d+\beta}}
\quad \text{if} ~\rho \geq (10)^{-1/\alpha} \\ &N \cdot
(10)^{(d+\beta)/\alpha}e^{-(d+\beta)((10)^{1/\alpha}\rho-1)} \quad
\text{if}~ \rho < (10)^{-1/\alpha}. \end{cases}
$$
Then, $\overline{\phi}_\beta$ is  continuously differentiable
  on $[0,\infty)$ such that
$$
\overline{\phi}_\beta(\infty)=0, \quad |\phi(x)| + |\nabla \phi(x)|
+ |x||\nabla \phi(x)| \leq \overline{\phi}_\beta(|x|), \quad
\int_0^\infty |\overline{\phi}_\beta'(\rho)|~d\rho \leq K
$$
and for each $r \geq (10)^{-1/\alpha}$,
$$
\int_r^\infty |\overline{\phi}_\beta'(\rho)| \rho^d~d\rho
=\int_r^\infty(d+\beta) \frac{N}{\rho^{d+1+\beta}} \rho^d~d\rho
=\frac{(d+\beta)N}{\beta} r^{-\beta}.
$$
The lemma is proved. \qed

\section{Some estimates on  $\cG f$}

In this section we develop some estimates of $\cG f$ by adopting the
approaches  in \cite{Kr01}, where the case $\alpha=2$ is studied.
Fix $f\in C_0^\infty(\bR^{d+1},H)$ and denote $u=\cG f$.

First, we prove a version of  Theorem \ref{main theorem} when $p=2$.

\begin{lemma}
\label{2-1}
 There exists a constant $N=N(\nu,\lambda,\alpha,K)>0$ so that for
 any $T\in (-\infty,\infty]$,
\begin{eqnarray}
\label{4023} \|u\|^2_{L_2(\bR^{d+1}\cap \{t\leq T\})} \leq
N\|f\|^2_{L_2(\bR^{d+1}\cap\{t \leq T\})}.
\end{eqnarray}
\end{lemma}

\pf  By the continuity of $f$, the range of $f$ belongs to a
separable subspace of $H$. Thus by using  a countable orthonormal
basis of this subspace and the Fourier transform one easily finds

\begin{eqnarray}
\notag \|u\|^2_{L_2(\bR^{d+1} \cap \{t \leq T\})}
&=&\int_{\bR^d} \int_{-\infty}^T [\int_{-\infty}^t |\hat{\psi}(\xi (t-s)^{1/\alpha})|^2|\hat{f}(s,\xi)|^2_H \frac{ds}{t-s}]dt d\xi\\
\notag &=&\int_{\bR^d} \int_{-\infty}^T \int_{-\infty}^T I_{s \leq t} |\hat{\psi}(\xi (t-s)^{1/\alpha})|^2|\hat{f}(s,\xi)|^2_H \frac{dt}{t-s}  ds d\xi\\
\label{411}&=&\int_{\bR^d} \int_{-\infty}^T \int_{0}^{T-s}
|\hat{\psi}(\xi t^{1/\alpha})|^2 \frac{dt}{t} |\hat{f}(s,\xi)|^2_H
ds d\xi.
\end{eqnarray}
By the assumption on $\psi$, for some $\nu, \lambda, K>0$,
$$
|\hat{\psi}(\xi)| \leq  K |\xi|^\nu,\quad
|\xi|^\lambda|\hat{\psi}(\xi)|\leq K.
$$
This and the change of  the variables $|\xi|^\alpha t \to t$ easily
lead to
\begin{eqnarray}
\int_0^\infty |\hat{\psi}(\xi t^{1/\alpha})|^2 \frac{dt}{t}
=\int_0^\infty |\hat{\psi}(t^{1/\alpha}\frac{\xi}{|\xi|})|^2
\frac{dt}{t}\nonumber\\
\leq K^2\int^1_0 t^{-1+2\nu/\alpha}dt+K^2\int^{\infty}_1
t^{-1-2\lambda/\alpha}dt \leq N(\nu, \alpha, \lambda,K).\label{412}
\end{eqnarray}
Plugging (\ref{412}) into (\ref{411}),
\begin{eqnarray*}
\|u\|^2_{L_2(\bR^{d+1} \cap \{t \leq T\})} \leq N \int_{-\infty}^T
\int_{\bR^d} |\hat{f}(s,\xi)|^2_H ~d\xi ds.
\end{eqnarray*}
The last expression is equal to the right-hand side of (\ref{4023}),
and therefore  the lemma is proved. \qed

For a real-valued function $h$ defined on $\bR^d$, define the
maximal function
$$
\bM_x h(x) := \sup_{r>0} \frac{1}{|B_r(x)|} \int_{B_r(x)} |h(y)| dy,
$$
where $| B_r(x) |$ denotes Lebesgue measure of $B_r(x)$. Similarly,
for measurable functions $h=h(t)$ on $\bR$ we introduce $\bM_th$ as
the maximal function of $h$ relative to symmetric intervals:
$$
\bM_t h(t) := \sup_{r>0} \frac{1}{2r} \int_{-r}^r |h(t+s)|\, ds.
$$
For a function $h(t,x)$ of two variables,  set
$$
\bM_x h(t,x) := \bM_x(h(t,\cdot))(x), \quad \bM_th(t,x) =
\bM_t(h(\cdot,x))(t).
$$

Denote
\begin{equation}
                                  \label{4024}
 Q_0 := [-2^\alpha,0] \times [-1,1]^d.
\end{equation}

\begin{corollary}
                                    \label{co1}
Assume that the support of $f$ is within $[-10,10] \times B_{3d}$. Then for any $(t,x) \in Q_0$
\begin{eqnarray}
\label{4025} \int_{Q_0} |u(s,y)|^2 \, dsdy \leq N  \bM_t \bM_x
|f|_H^2 (t,x),
\end{eqnarray}
where $N$ depends only on $d,\alpha,\nu,\lambda$ and $K$.
\end{corollary}

\pf By the Lemma \ref{2-1},
$$
\int_{Q_0} |u(s,y)|^2 ~dsdy \leq \int_{-\infty}^0 \int_{\bR^d}
|u(s,y)|^2  dyds \leq N \int_{-10}^0 \int_{B_{3d}} |f(s,y)|_H^2
dyds.
$$
Since  $|x-y| \leq |x|+|y| \leq   4d$ for any $(t,x) \in Q_0$ and $y
\in B_{3d}$,
\begin{eqnarray*}
 \int_{-10}^0 \int_{B_{3d}} |f(s,y)|_H^2 dyds
\leq  \int_{-10}^0 \int_{|x-y| \leq 4d} |f(s,y)|_H^2 dyds
&\leq& N \int_{-10}^0 \bM_x|f(s,x)|_H^2 ds \\
&\leq& N \bM_t\bM_x f(t,x).
\end{eqnarray*}
The lemma is proved. \qed

We generalize Corollary \ref{co1} as follows.

 \begin{lemma}
                                        \label{2-3}
 Assume that $f(t,x)=0$ for $t \neq (-10,10)$. Then for any $(t,x) \in
 Q_0$,
\begin{eqnarray*}
 \int_{Q_0} |u(s,y)|^2 ~dsdy \leq N \bM_t \bM_x |f|_H^2 (t,x),
\end{eqnarray*}
 where $N=N(d,\alpha,\nu,\lambda,\delta,K)$.
\end{lemma}

\pf First, notice that if  $0 \leq \varepsilon \leq R \leq \infty$,
and $F$ and $G$ are smooth enough, then
\begin{eqnarray}
\notag \int_{R \geq |z| \geq \varepsilon} F(z) G(|z|)~dz
= - \int_\varepsilon^R G'(\rho)(\int_{|z| \leq \rho} F(z) dz) d\rho \\
+ G(R) \int_{|z| \leq R} F(z)dz - G(\varepsilon) \int_{|z| \leq
\varepsilon} F(z)dz. \label{4026}
\end{eqnarray}
Indeed, (\ref{4026}) is obtained by applying integration by parts to
$$
\int_\varepsilon^R G(\rho)\frac{d}{d\rho} \left(\int_{B_\rho(0)}
F(z)\, dz \right)d\rho=\int_\varepsilon^R G(\rho)
\left(\int_{\partial B_\rho(0)} F(s) \, dS_{\rho} \right) d\rho
=\int_{R \geq |z| \geq \varepsilon} F(z) G(|z|)\,dz.
$$

   Now take  $\zeta \in C_0^\infty (\bR^d)$ such that $\zeta =1$ in
$B_{2d}$ and $\zeta=0$ outside of $B_{3d}$. Set $\cA = \zeta f$ and
$\cB = (1-\zeta)f$. By Minkowski's inequality, $\cG f \leq \cG \cA+
\cG \cB$. Since $\cG \cA$ can be estimated by Corollary \ref{co1},
we may assume that $f(t,x)=0$ for $x \in B_{2d}$.

Denote $\overline{f}=|f|_H$, take $0 > s > r > -10$, and see
\begin{eqnarray*}
|\Psi_{s-r}f(r, \cdot)(y)|_H
&\leq&  (s-r)^{-d/\alpha}\int_{\bR^d} |\psi(z/(s-r)^{1/\alpha})| |f(r,y-z)|_H ~dz\\
&\leq&  (s-r)^{-d/\alpha}\int_{\bR^d}
\overline{\psi}(|z|/(s-r)^{1/\alpha}) \overline{f}(r,y-z)~dz.
\end{eqnarray*}
Observe that if $(s,y) \in Q_0$ and $|z| \leq \rho$ with a $\rho >
1$, then
\begin{eqnarray}
\label{4027} |x-y| \leq 2d, \quad B_\rho(y) \subset B_{2d+\rho}(x)
\subset B_{\mu \rho}(x), \quad \mu=2d+1,
\end{eqnarray}
whereas if $|z| \leq 1$, then $|y-z| \leq  2d$ and $f(r,y-z)=0$.
Thus by (\ref{4026}), for $0 >s > r > -10$ and $(s,y) \in Q_0$
\begin{eqnarray*}
|\Psi_{s-r}f(r, \cdot)(y)|_H
&\leq&  (s-r)^{-(d+1)/\alpha}\int_1^\infty |\overline{\psi}'(\rho/(s-r)^{1/\alpha})|(\int_{|z| \leq \rho} \overline{f}(r,y-z)~dz) ~d\rho\\
&=&  (s-r)^{-(d+1)/\alpha}\int_1^\infty |\overline{\psi}'(\rho/(s-r)^{1/\alpha})|(\int_{B_{\rho}(y)} \overline{f}(r,z)~dz) ~d\rho\\
&\leq&  (s-r)^{-(d+1)/\alpha}\int_1^\infty |\overline{\psi}'(\rho/(s-r)^{1/\alpha})|(\int_{B_{\mu \rho}(x)} \overline{f}(r,z)~dz) ~d\rho\\
&\leq&  N \bM_x\overline{f}(r,x)(s-r)^{-(d+1)/\alpha}\int_1^\infty |\overline{\psi}'(\rho/(s-r)^{1/\alpha})|\rho^d~ d\rho\\
&=&  N  \bM_x\overline{f}(r,x)\int_{(s-r)^{-1/\alpha}}^\infty
|\overline{\psi}'(\rho)|\rho^d~ d\rho\leq  N
\bM_x\overline{f}(r,x)(s-r)^{\delta/\alpha},
\end{eqnarray*}
where the last inequality follows from (\ref{eqn 06.08}) and the
inequality $(s-r)^{-1/\alpha}\geq 10^{-1/\alpha}$. By Jensen's
inequality $(\bM_x\overline{f})^2 \leq \bM_x \overline{f}^2$, and
therefore, for  any $(s,y) \in Q_0$ (remember $\delta\geq \alpha/2$)
\begin{eqnarray*}
|u(s,y)|^2
= \int_{-\infty}^s |\Psi_{s-r}f(r,\cdot)(y)|_H^2 \frac{dr}{s-r}
&\leq& N \int_{-10}^s \bM_x \overline{f}^2(r,x) (s-r)^{2\delta/\alpha -1}dr\\
&\leq& N \int_{-10}^0 \bM_x \overline{f}^2(r,x) dr
\leq  N \bM_t\bM_x \overline{f}^2(t,x).
\end{eqnarray*}
The lemma is proved.\qed
\begin{lemma}
\label{2-4}
 Assume that $f(t,x)=0$ for $t \geq -8$. Then for any $(t,x) \in Q_0$
\begin{eqnarray}
                                                 \label{4028}
\int_{Q_0} |u(s,y) -u(t,x)|^2~dsdy \leq N\bM_t \bM_x |f|_H^2 (t,x),
\end{eqnarray}
where $N=N(d,\alpha,\nu,\lambda,\delta,K)$.
\end{lemma}

\pf Obviously it is enough to show that
\begin{equation}
                    \label{eqn 6.08.8}
 \sup_{Q_0} [ |D_su|^2 + |\nabla u|^2] \leq N \bM_t \bM_x
|f|_H^2 (t,x).
\end{equation}
By Minkowski's inequality  the derivative of a norm is less than or
equal to the norm of the derivative  if both exist. Thus for fixed
$(s,y) \in Q_0$ we have
$$
|\nabla u (s,y)|^2 \leq \int_{-\infty}^{-8} |\nabla
\Psi_{s-r}f(r,\cdot)(y)|^2_H
\frac{dr}{s-r}=:\int^{-8}_{-\infty}I^2(r,s,y)\frac{dr}{s-r},
$$
where
\begin{eqnarray*}
I(r,s,y)&:=& |\nabla \Psi_{s-r}f(r,\cdot)(y)|_H\\
&=&(s-r)^{-(d+1)/\alpha}|\int_{\bR^d}(\nabla \psi)(z/(s-r)^{1/\alpha})f(r,y-z)~dz|_H\\
&\leq& (s-r)^{-(d+1)/\alpha}\int_{\bR^d}
\overline{\psi}(|z|/(s-r)^{1/\alpha})\overline{f}(r,y-z)~dz=:\tilde{I}(r,s,y),
\end{eqnarray*}
and $\bar{f}:=|f|_{H}$.  Using (\ref{4026}) and (\ref{4027}) again,
we get
 for $s>r$,
\begin{eqnarray*}
\tilde{I}(r,s,y) &\leq&(s-r)^{-(d+2)/\alpha}\int_0^\infty
\overline{\psi}'(\rho/(s-r)^{1/\alpha})
(\int_{B_\rho(y)}\overline{f}(r,z)~dz)~d\rho\\
&\leq&(s-r)^{-(d+2)/\alpha}\int_0^\infty \overline{\psi}'(\rho/(s-r)^{1/\alpha})
(\int_{B_{2d+\rho}(x)}\overline{f}(r,z)~dz)~d\rho\\
&\leq& N \bM_x\overline{f}(r,x) (s-r)^{-(d+2)/\alpha}\int_0^\infty
\overline{\psi}'(\rho/(s-r)^{1/\alpha})
(2d + \rho)^d d\rho\\
&=& N  \bM_x\overline{f}(r,x) (s-r)^{-1/\alpha}\int_0^\infty
\overline{\psi}'(\rho)(2d/(s-r)^{1/\alpha} + \rho)^d d\rho.
\end{eqnarray*}
For $r \leq -8$, we have $s-r\geq 2^\alpha $ and
\begin{eqnarray*}
\int_0^\infty |\overline{\psi}'(\rho)| ( 2d / (s-r)^{1/\alpha} +
\rho)^d~d\rho \leq \int_0^\infty |\overline{\psi}'(\rho)| ( d+
\rho)^d~d\rho \leq N,
\end{eqnarray*}
\begin{eqnarray*}
\tilde{I}(r,s,y) \leq N \bM_x\overline{f}(r,x) (s-r)^{-1/\alpha}
\end{eqnarray*}
and
\begin{eqnarray*}
|\nabla u(s,y)|^2 \leq \int_{-\infty}^{-8} \tilde{I}^2(r,s,y)
\frac{dr}{s-r}
&\leq& N \int_{-\infty}^{-8} \bM_x\overline{f}^2(r,x)  \frac{dr}{(s-r)^{2/\alpha+1}}\\
&\leq& N \int_{-\infty}^{-8} \bM_x\overline{f}^2(r,x)
\frac{dr}{(-4-r)^{2/\alpha+1}}.
\end{eqnarray*}
By the integration by parts,
\begin{eqnarray}
|\nabla u(s,y)|^2 &\leq& N\int_{-\infty}^{-8} \tilde{I}^2(r,s,y)
\frac{dr}{s-r} \nonumber\\
&\leq& N\int_{-\infty}^{-8} \frac{1}{(-4-r)^{2/\alpha+2}} (\int_r^0 \bM_x \overline{f}^2(p,x)~dp)~dr \nonumber\\
&\leq& N \bM_t\bM_x\overline{f}^2(t,x) \int_{-\infty}^{-8}
\frac{|r|}{(-4-r)^{2/\alpha+2}} ~dr=N\bM_t\bM_x\overline{f}^2(t,x)
\label{eqn 06.08.2}.
\end{eqnarray}

To estimate $D_su$, we proceed similarly. By Minkowski's inequality,
\begin{eqnarray}
|D_su(s,y)|^2 &\leq& N \int_{-\infty}^{-8} \left(
|D_s\Psi_{s-r}f(r,y)|_H^2
\frac{1}{s-r}+ |\Psi_{s-r}f(r,y)|_H^2\frac{1}{(s-r)^3}\right)\,dr \nonumber \\
&=:&  N \int_{-\infty}^{-8} J^2(r,s,y) \frac{1}{s-r}\,dr, \label{eqn
6.08.5}
\end{eqnarray}
where
\begin{eqnarray*}
J(r,s,y)&:=&(s-r)^{-d/\alpha} |\int_{\bR^d}D_s \psi(z/(s-r)^{1/\alpha}) f(r,y-z)~dz|_H \\
&&+(s-r)^{-d/\alpha-1} | \int_{\bR^d} \psi(z/(s-r)^{1/\alpha}) f(r,y-z)~dz|_H\\
&=&(s-r)^{-d/\alpha} |\int_{\bR^d} \nabla \psi(z/(s-r)^{1/\alpha}) \cdot (-\frac{1}{\alpha}(s-r)^{-1/\alpha-1}z ) f(r,y-z)~dz|_H\\
&&+(s-r)^{-d/\alpha-1} | \int_{\bR^d} \psi(z/(s-r)^{1/\alpha}) f(r,y-z)~dz|_H\\
&\leq& N (s-r)^{-d/\alpha-1 } \int_{\bR^d}
\overline{\psi}(|z|/(s-r)^{1/\alpha})
\overline{f}(r,y-z)~dz=N\tilde{I}(r,s,y).
\end{eqnarray*}
This, (\ref{eqn 06.08.2}) and  (\ref{4028}) lead to (\ref{eqn
6.08.8}). The lemma is proved. \qed

\mysection{Proof of Theorem \ref{main theorem}}
                   \label{proof of theorem}
Note that we may assume $a=-\infty$ and $b=\infty$. Indeed, for any
$f\in C^{\infty}_0((a,b)\times \bR^d,H)$ we have $f\in
C^{\infty}_0(\bR^{d+1},H)$, and inequality (\ref{4022}) with
$a=-\infty$ and $b=\infty$ implies the inequality with any pair of
$(a,b)$. Since in this case the theorem is already proved if $p=2$,
we assume $p>2$.

Let $\cF$ be the collections of all balls $Q\subset \bR^{d+1}$ of
the type
$$
Q_{c}(s,y):=\{ (s-c^\alpha,s) \times (y^1-c/2 , y^1 +c/2) \cdots
(y^d-c/2, y^d +c/2) \}, \quad c>0.
$$
For a measurable function $h(t,x)$  on $\bR^{d+1}$, define the sharp
function
$$
h^\#(t,x) := \sup_{Q} \frac{1}{|Q|} \int_Q |h (t,x)-h_Q|\, dyds,
$$
where
$$
h_Q = \dashint_Q  h ~dyds:=\frac{1}{|Q|} \int_Q h(s,y)\, dyds
$$
and the supremum is taken over all balls $Q\in \cF$  containing
$(t,x)$.

\begin{thm}(Fefferman-Stein). For any $1<q<\infty$ and $h\in L_q(\bR^{d+1})$,
\begin{equation}
                          \label{eqn 6.08.8}
\|h\|_{L^q}\leq N(q)\|h^\#\|_{L^q}.
\end{equation}
\end{thm}
\pf Inequality (\ref{eqn 6.08.8}) is a consequence of Theorem IV.2.2
in \cite{Ste}, because  the balls $Q_c(s,y)$ satisfy the
conditions (i)-(iv) in section 1.1 of \cite{Ste} :\\
(i) $Q_{c}(t,x)\cap Q_{c}(s,y)\neq \emptyset$ implies
$Q_{c}(s,y)\subset Q_{N_1c}(t,x)$ ;\\
(ii) $|Q_{N_1c}(t,x)|\leq N_2 |Q_{c}(t,x)|$ ;\\
(iii) $\cap_{c>0}\overline{Q}_{c}(t,x)=\{(t,x)\}$ and
$\cup_{c}Q_{c}(t,x)=\bR^{d+1}$ ;\\
(iv) for each open set $U$ and $c>0$, the function $(t,x)\to
|Q_{c}(t,x)\cap U|$ is continuous.  \qed

Next we prove
\begin{equation}
                        \label{eqn 6.08.9}
  (\cG f)^{\#}(t,x)\leq N (\bM_t\bM_x|f|^2_{H})^{1/2}(t,x).
  \end{equation}
By Jensen's inequality, to prove (\ref{eqn 6.08.9}) it suffices to
prove that for each $Q=Q_{c}(s,y)\in \cF$ and $(t,x)\in Q$,
\begin{eqnarray}
\label{4033} \dashint_Q |\cG f- (\cG f)_Q|^2~dyds \leq
N(d,\alpha,\nu,\lambda,\delta,K) \bM_t\bM_x|f|_H^2(t,x).
\end{eqnarray}
It is easy to check that to prove (\ref{4033}) we may assume
$(s,y)=(0,0)$. Note that for any $c
> 0$, $\Psi_th(c~ \cdot)(x) = \Psi_{tc^\alpha}h(cx)$ and
\begin{eqnarray}
\notag \cG f(c^{\alpha}~ \cdot , c ~\cdot)(t,x)
\notag  &=&[\int_{-\infty}^t |\Psi_{(t-s)c^\alpha} f(c^{\alpha} s,\cdot)(cx)|_H^2 \frac{ds}{t-s}]^{1/2}\\
\notag  &=&[\int_{-\infty}^{t c^{\alpha}} |\Psi_{(t-c^{-\alpha}s)c^\alpha} f(s,\cdot)(cx)|_H^2 \frac{ c^{-\alpha} ds}{t-c^{-\alpha}s}]^{1/2}\\
\notag  &=&[\int_{-\infty}^{t c^{\alpha}} |\Psi_{(c^\alpha t-s)} f(s,\cdot)(cx)|_H^2 \frac{ds}{c^\alpha t-s}]^{1/2}\\
\label{5271}&=&\cG f(c^\alpha t, cx).
\end{eqnarray}
Since dilations don't affect averages, (\ref{5271}) shows that it
suffices to prove (\ref{4033}) when $c=2$, that is $Q=Q_0$ from
(\ref{4024}).  Now we take a function $\zeta \in C_0^\infty(\bR)$
such that $\zeta=1$ on $[-8,8]$, $\zeta=0$ outside of $[-10,10]$,
and $1 \geq \zeta \geq 0$. Define
$$
\cA(s,y) := f(s,y) \zeta(s), \quad \cB(s,y):=
f(s,y)-\cA(s,y)=f(s,y)(1-\zeta(s)).
$$
Then
$$
\Psi_{t-s}\cA(s,\cdot) = \zeta(s) \Psi_{t-s}f(s,\cdot), \ \quad \cG
f \leq \cG \cA + \cG \cB,  \quad \cG \cB \leq \cG f
$$ and  for
any constant $c$, $|\cG f -c | \leq |\cG\cA| + |\cG \cB -c|$. Thus
\begin{eqnarray*}
\dashint_{Q_0} |\cG f- (\cG f)_{Q_0}|^2~dyds \leq 4 \dashint_{Q_0}
|\cG f-c|^2~dyds \leq 8 \dashint_{Q_0} |\cG \cA|^2~dyds +
8\dashint_{Q_0} |\cG \cB - c |^2~dyds.
\end{eqnarray*}
Taking $c=\cG \cB(t,x)$, from Lemma \ref{2-3}  we get
\begin{eqnarray*}
 \dashint_{Q_0} |\cG f- (\cG f)_{Q_0}|^2~dyds
&\leq& 8 \dashint_{Q_0} |\cG \cA|^2~dyds + 8\dashint_{Q_0} |\cG \cB - \cG \cB (t,x) |^2~dyds\\
&\leq& N \bM_t \bM_x |f|_H^2 (t,x)+8\dashint_{Q_0} |\cG \cB - \cG
\cB (t,x) |^2~dyds.
\end{eqnarray*}
In addition, setting $f_1(s, y) := \cB(s,y)$ on $s  \leq 0$ and $f_1(s,y):=0$ on $s >0$, from Lemma \ref{2-4} we see
\begin{eqnarray*}
 \bM_t \bM_x |f|_H^2 (t,x)+8\dashint_{Q_0} |\cG \cB - \cG \cB (t,x) |^2~dyds
&\leq&  \bM_t \bM_x |f|_H^2 (t,x)+N \bM_t \bM_x |f_1|_H^2 (t,x)\\
&\leq&  N\bM_t \bM_x |f|_H^2 (t,x).
\end{eqnarray*}
This proves (\ref{eqn 6.08.9}).

Finally, combining the Fefferman-Stein theorem and Hardy-Littlewood
maximal theorem (see, for instance, \cite{Ste}), we conclude (recall
$p/2
>1$)
\begin{eqnarray*}
\|u\|_{L_p(\bR^{d+1})}^p \leq N \|(\bM_t\bM_x|f|_H^2)^{1/2}\|_{L_p(\bR^{d+1})}^p
&=&N \int_{\bR^d} \int_\bR (\bM_t \bM_x |f|_H^2)^{p/2}~dt~dx \\
&\leq& N \int_{\bR^d}  \int_\bR (\bM_x |f|_H^2 )^{p/2}~dt~dx \\
& =& N \int_{\bR}  \int_{\bR^d} (\bM_x |f|_H^2 )^{p/2}~dx~dt \\
&\leq&  N \|f\|_{L_p(\bR^{d+1},H)}^p.
\end{eqnarray*}
The theorem is proved. \qed

Below we explain why Theorem \ref{main theorem} implies (\ref{for
fun}). By (\ref{eqn 333}) and Remark \ref{remark 22}, for any
solution $u$ of (\ref{eqn 0}), we have
\begin{equation}
                         \label{eqn 6.9.11}
\bE \int^T_0\|\partial^{\alpha/2}_xu(t,\cdot)\|^p_{L_p}dt \leq N \bE
\int_{\bR^d} \int_0^T |f|^p_{\ell_2}~dtdx.
\end{equation}
By (\ref{eqn 222}) and Burkholder-Davis-Gundy inequality,
\begin{equation}
                    \label{eqn 33333}
\bE \int^T_0\|u(t,\cdot)\|^p_{L_p}dt \leq N \bE
\int^T_0\int_{\bR^d}\left[\int^t_0|T_{\alpha,t-s}f(s,\cdot)(x)|^2_{\ell^2}
ds\right]^{p/2}dxdt,
\end{equation}
 and  by Jensen's
inequality
\begin{equation}
                         \label{eqn 6.9.10}
|T_{t-s}f(x)|^2_{\ell_2}=\sum_k\left(\int_{\bR^d}p(t-s,y)f^k(x-y)dy\right)^2\leq
N (p(t-s,\cdot)*|f(\cdot)|^2_{\ell_2})(x).
\end{equation}
Thus (\ref{eqn 6.9.11}), (\ref{eqn 33333}), (\ref{eqn 6.9.10}) and
Remark \ref{remark 6.14} imply
$$
\bE\int^T_0\|u\|^p_{H^{\alpha/2}_p}dt\leq N
\bE\int^T_0(\|u\|^p_{L^p}+\|\partial^{\alpha/2}u\|^p_{L^p})dt\leq
N(T,d,\alpha)\bE \int_{\bR^d} \int_0^T |f|^p_{\ell_2}~dtdx.
$$

\end{document}